\documentclass[12pt]{amsart}

\usepackage{amsfonts}
\usepackage{amsthm}
\usepackage[psamsfonts]{amssymb}
\usepackage[dvipdfmx]{graphicx}
\usepackage{ascmac}
\usepackage{amsmath}
\usepackage{fancybox}
\usepackage{fancyhdr}
\usepackage{enumerate}
\usepackage{calrsfs}
\usepackage{mathrsfs}
\usepackage{color}
\usepackage{multicol}
\usepackage{framed}

\theoremstyle{plain}
\newtheorem{thm}[equation]{Theorem}

\newtheorem{defi}[equation]{\it Definition}
\newtheorem{cor}[equation]{Corollary}
\newtheorem{rem}[equation]{\it Remark}

\newtheorem*{Acknowledgements}{\it Acknowledgements}

\address{Department of Mathematics, Graduate School of Science, Osaka University, 1-1 Machikaneyama, Toyonaka, Osaka 560-0043, JAPAN}
\email{h-tadano@cr.math.sci.osaka-u.ac.jp}
\title{A lower diameter bound for compact domain manifolds of shrinking Ricci-harmonic solitons}
\dedicatory{Dedicated to Professor Akito Futaki on the occasion of his Kanreki (60th birthday)}
\author{Homare TADANO}
\date{submitted June 12, 2014, revised October 14, 2014}
\thanks{This work was supported by Moriyasu Graduate Student Scholarship Foundation}

\begin{document}

\begin{abstract}
In this paper, we shall give a lower diameter bound for compact domain manifolds of shrinking Ricci-harmonic solitons. Our result may be regarded as a generalization to Ricci-harmonic geometry of the recent works by Fern\'{a}ndez-L\'{o}pez and Garc\'{i}a-R\'{i}o (Q. J. Math. 61, 319--327, 2010), Futaki and Sano (Asian J. Math. 17, 17--32, 2013), and Futaki \textit{et al} (Ann. Global Anal. Geom. 44, 105--114, 2013).
\end{abstract}

\maketitle

\numberwithin{equation}{section}

\section{Introduction}

Recently, various geometric flows have been studied extensively. In this paper, we shall deal with a generalization of the Ricci flow.

\subsection{The Bernhard List's flow}

Let $(M, g(t))$ be a family of $m$-dimensional Riemannian manifolds with Riemannian metrics $g(t)$ evolving by the following coupled system:
\begin{equation}\label{List-Flow}
\left\{
\begin{aligned}
\frac{\partial}{\partial t} g(t) & = - 2 \operatorname{Ric}_{g(t)} + 2 \alpha d \phi(t) \otimes d \phi(t), \\
\frac{\partial}{\partial t} \phi(t) & = \Delta_{g(t)} \phi(t), 
\end{aligned}
\right.
\end{equation}
where $\operatorname{Ric}_{g(t)}$ denotes the Ricci tensor with respect to the evolving metric $g(t)$, $\alpha \geqslant 0$ is a non-negative constant, $\phi(t) : (M, g(t)) \rightarrow \mathbb{R}$ is a family of smooth functions on $M$, and $\Delta_{g(t)} := g^{ij}(t) \nabla_{i} (t) \nabla_{j}(t)$ denotes the Laplace-Beltrami operator with respect to $g(t)$. The flow (\ref{List-Flow}) is called a \textit{Bernhard List's flow} and was introduced by List \cite{List1, List2}. The short time existence is proved. A typical example would be the Ricci flow \cite{Hamilton} playing an important role in the Perelman's work \cite{Perelman}, in which case $\phi(t) : (M, g(t)) \rightarrow \mathbb{R}$ is a constant function. The motivation to study the Bernhard List's flow stems from its connection to general relativity. The stationary points of the flow correspond to the static Einstein vacuum equations \cite{List1, List2}.

\subsection{The Ricci-harmonic flow}

After List introduced the Bernhard List's flow, a new geometric flow was introduced. Let $(M, g(t))$ be a family of $m$-dimensional Riemannian manifolds with Riemannian metrics $g(t)$ evolving by the following coupled system:
\begin{equation}\label{Ricci-Harmonic-Flow}
\left\{
\begin{aligned}
\frac{\partial}{\partial t} g(t) & = - 2 \operatorname{Ric}_{g(t)} + 2 \alpha(t) \nabla \phi(t) \otimes \nabla \phi(t), \\
\frac{\partial}{\partial t} \phi(t) & = \tau_{g(t)} \phi(t), 
\end{aligned}
\right.
\end{equation}
where $\alpha(t) \geqslant 0$ is a non-negative time-dependent constant, $\phi(t) : (M, g(t)) \rightarrow (N, h)$ is a family of smooth maps between $(M, g(t))$ and a fixed Riemannian manifold $(N, h)$ of dimension $n$, $\nabla \phi(t) \otimes \nabla \phi(t) := \phi(t)^{*} h$ is the pull-back of the metric $h$ via $\phi(t)$, and $\tau_{g(t)} \phi(t) := \mathrm{trace} \nabla d \phi(t)$ denotes the tension field of $\phi(t)$ with respect to $g(t)$. The flow (\ref{Ricci-Harmonic-Flow}) is called a \textit{Ricci-harmonic flow} and was introduced by M\"{u}ller \cite{Muller1, Muller3}. The short time existence is proved. The Bernhard List's flow (\ref{List-Flow}) is an example of the Ricci-harmonic flow, in which case $(N, h) = (\mathbb{R}, dr^{2})$. More examples can be found in \cite{Muller1, Muller3}. By denoting $S(t) := \operatorname{Ric}_{g(t)} - \alpha d \phi(t) \otimes d \phi(t)$ and $S(t) := \operatorname{Ric}_{g(t)} - \alpha(t) \nabla \phi(t) \otimes \nabla \phi(t)$ in (\ref{List-Flow}) and (\ref{Ricci-Harmonic-Flow}) respectively, the first equations in (\ref{List-Flow}) and (\ref{Ricci-Harmonic-Flow}) are written as
\begin{equation}\label{Flow-S}
\frac{\partial}{\partial t} g(t) = -2 S(t).
\end{equation}
As with the Ricci flow, under (\ref{List-Flow}) and (\ref{Ricci-Harmonic-Flow}), some differential Harnack inequalities for several heat-type equations are obtained, respectively \cite{Bailesteanu, Fang1, Fang3, Zhao, Zhu}. Note that the papers \cite{Fang2, Guo-He, Guo-Ishida1, Guo-Ishida2, Ishida} study these Harnack inequalities under more general settings, namely, under the flow (\ref{Flow-S}) for smooth symmetric two-tensors $S(t)$ with some technical assumptions on evolving tensor quantities associated to $S(t)$.

\subsection{The Ricci-harmonic soliton}

Let $(M, g)$ and $(N, h)$ be two static Riemannian manifolds of dimension $m$ and $n$, respectively and let $\phi : (M, g) \rightarrow (N, h)$ be a smooth map between the domain manifold $(M, g)$ and the target manifold $(N, h)$, $f : M \rightarrow \mathbb{R}$ a $\mathcal{C}^{2}$-class function on $M$, and $\lambda \in \mathbb{R}$ a real number.

\begin{defi}[Williams \cite{Williams}]\rm
The $4$-tuple $((M, g), (N, h), \phi, \lambda)$ is called \textit{harmonic-Einstein} if it satisfies the following coupled system:
\[
\left\{
\begin{aligned}
& \operatorname{Ric}_{g} - \alpha \nabla \phi \otimes \nabla \phi = \lambda g, \\
& \tau_{g} \phi = 0, 
\end{aligned}
\right.
\]
where $\operatorname{Ric}_{g}$ denotes the Ricci tensor with respect to the metric $g$, $\alpha \in \mathbb{R}$ is a constant, $\nabla \phi \otimes \nabla \phi := \phi^{*} h$ is the pull-back of the metric $h$ via $\phi$, and $\tau_{g} \phi := \mathrm{trace} \nabla d \phi$ denotes the tension field of $\phi$ with respect to $g$.
\end{defi}

\begin{defi}[M\"{u}ller \cite{Muller1, Muller3}]\rm
The $5$-tuple $((M, g), (N, h), \phi, f, \lambda)$ is called a \textit{Ricci-harmonic soliton} if it satisfies the following coupled system:
\begin{equation}\label{Ricci-Harmonic-Soliton}
\left\{
\begin{aligned}
& \operatorname{Ric}_{g} - \alpha \nabla \phi \otimes \nabla \phi + \operatorname{Hess} f = \lambda g, \\
& \tau_{g} \phi = \left< \nabla \phi, \nabla f \right>, 
\end{aligned}
\right.
\end{equation}
where $\alpha \geqslant 0$ is a non-negative constant and $\operatorname{Hess} f$ denotes the Hessian of $f$. We say that the soliton $((M, g), (N, h), \phi, f, \lambda)$ is \textit{shrinking}, \textit{steady}, and \textit{expanding} described as $\lambda > 0, \lambda = 0$, and $\lambda < 0$, respectively. If $f$ is constant in (\ref{Ricci-Harmonic-Soliton}), then the soliton is harmonic-Einstein. In such a case, we say that the soliton is \textit{trivial}.
\end{defi}

The soliton (\ref{Ricci-Harmonic-Soliton}) above is a self-similar solution for the coupled system (\ref{Ricci-Harmonic-Flow}). Note that, if $(N, h) = (\mathbb{R}, dr^{2})$ and $\phi : (M, g) \rightarrow (\mathbb{R}, dr^{2})$ is a constant function in (\ref{Ricci-Harmonic-Soliton}), then the soliton is exactly a gradient Ricci soliton. Since the Ricci-harmonic flow is a generalization of the Ricci flow, a natural question to ask is whether fundamental theorems as in the Ricci flow also hold for the Ricci-harmonic flow. In this direction, corresponding theories have been established, such as, the extension theorem \cite{Cheng-Zhu}, the no breathers theorem, the non-collapsing theorem \cite{Muller1, Muller3}, the monotonicity formula \cite{Li, Muller2}, and the volume growth estimate \cite{Yang-Shen}. As with the Ricci soliton, any non-trivial Ricci-harmonic soliton $((M, g), (N, h), \phi, f, \lambda)$ with a compact domain manifold $M$ is shrinking \cite{Williams}.

\medskip

In this paper, we give a lower diameter bound for compact domain manifolds of shrinking Ricci-harmonic solitons. Our main result is the following:

\begin{thm}\label{Main-Theorem}
Let $((M, g), (N, h), \phi, f, \lambda)$ be a non-trivial shrinking Ricci-harmonic soliton satisfying {\rm (\ref{Ricci-Harmonic-Soliton})}. Suppose that the domain manifold $M$ is compact. Then the diameter of the domain manifold $(M, g)$ has the universal lower bound
\begin{equation}\label{Main-Theorem-Ineq}
\operatorname{diam}(M, g) \geqslant \frac{2(\sqrt{2} - 1)\pi}{\sqrt{\lambda}}.
\end{equation}
\end{thm}

\begin{cor}
Let $((M, g), (N, h), \phi, f, \lambda)$ be a shrinking Ricci-harmonic soliton satisfying {\rm (\ref{Ricci-Harmonic-Soliton})}. Suppose that the domain manifold $M$ is compact. If the diameter of the domain manifold $(M, g)$ satisfies
\[
\operatorname{diam}(M, g) < \frac{2(\sqrt{2} - 1)\pi}{\sqrt{\lambda}}, 
\]
then the soliton must be harmonic-Einstein.
\end{cor}

\begin{rem}\rm
In Theorem \ref{Main-Theorem} above, if $(N, h) = (\mathbb{R}, dr^{2})$ and $\phi : (M, g) \rightarrow (\mathbb{R}, dr^{2})$ is a constant function, then the soliton $((M, g), (N, h), \phi, f, \lambda)$ appears as a compact shrinking gradient Ricci soliton and (\ref{Main-Theorem-Ineq}) recovers the lower diameter bound for compact shrinking Ricci solitons given by Futaki \textit{et al} \cite{Futaki-Li-Li}.
\end{rem}

\begin{Acknowledgements}\rm
I would like to thank Professors Toshiki Mabuchi and Kimio Miyajima for their encouragements. I also thank Professor Masashi Ishida for his comments.
\end{Acknowledgements}

\section{Proof of Theorem \ref{Main-Theorem}}\label{Section-2}

The proof of Theorem \ref{Main-Theorem} is almost the same as \cite{Futaki-Sano, Futaki-Li-Li}. By replacing the Ricci tensor $\operatorname{Ric}_{g}$ with the symmetric two-tensor $\operatorname{Ric}_{g} - \alpha \nabla \phi \otimes \nabla \phi$ in the argument of \cite{Futaki-Sano, Futaki-Li-Li}, we can give a proof of Theorem \ref{Main-Theorem} quite parallelly as the case of compact shrinking Ricci solitons.

\begin{proof}
The Witten-Laplacian $\Delta_{f}$ acting on $\mathcal{C}^{2}$-class functions on $(M, g)$ is defined by
\[
\Delta_{f} := \Delta - \nabla f \cdot \nabla, 
\]
where $\Delta = g^{ij} \nabla_{i} \nabla_{j}$ denotes the Laplace-Beltrami operator with respect to the metric $g$. We first show that $2 \lambda$ is an eigenvalue of the Witten-Laplacian $\Delta_{f}$. By taking the trace of the first equation in (\ref{Ricci-Harmonic-Soliton}), we have
\begin{equation}\label{Fundamental-Equation-1}
R - \alpha | \nabla \phi |^{2} + \Delta f = n \lambda, 
\end{equation}
where $R := g^{ij} R_{ij}$ denotes the scalar curvature on $(M, g)$. On the other hand, by taking a covariant derivative of the first equation in (\ref{Ricci-Harmonic-Soliton}), we obtain
\[
\nabla_{k} R_{ij} - \alpha \nabla_{k} (\nabla_{i} \phi \nabla_{j} \phi) + \nabla_{k} \nabla_{i} \nabla_{j} f = 0.
\]
By subtracting the same equation with indices $i$ and $k$ interchanged, we have
\[
\nabla_{k} R_{ij} - \nabla_{i} R_{kj} - \alpha (\nabla_{i} \phi \nabla_{k} \nabla_{j} \phi - \nabla_{k} \phi \nabla_{i} \nabla_{j} \phi) + R_{kijp} \nabla_{p} f = 0.
\]
Tracing just above with $g^{kj}$ yields
\[
\nabla_{j} R_{ij} - \nabla_{i} R - \alpha \left( \nabla_{i} \phi \tau_{g} \phi - \frac{1}{2} \nabla_{i} | \nabla \phi |^{2} \right) + R_{ip} \nabla_{p} f = 0.
\]
By using the contracted second Bianchi identity $\nabla_{j} R_{ij} = \frac{1}{2} \nabla_{i} R$ and plugging in both equations of (\ref{Ricci-Harmonic-Soliton}) for $R_{ip}$ and for $\tau_{g} \phi$, we obtain
\[
- \frac{1}{2} \nabla_{i} (R - \alpha | \nabla \phi |^{2} + | \nabla f |^{2} - 2 \lambda f) = 0.
\]
Hence, there exists a real constant $K \in \mathbb{R}$ such that
\begin{equation}\label{Fundamental-Equation-2}
R - \alpha | \nabla \phi |^{2} + | \nabla f |^{2} - 2 \lambda f = K.
\end{equation}
By combining two equalities (\ref{Fundamental-Equation-1}) and (\ref{Fundamental-Equation-2}), we have
\begin{equation}\label{Eigenvalue-1}
\Delta_{f} f = \Delta f - | \nabla f |^{2} = - 2 \lambda f + K', 
\end{equation}
where $K' := n \lambda - K$. By adding some constant on $f$, we may normalize $f$ such that
\[
\int_{M} f e^{- f} d \mathrm{vol}_{g} = 0.
\]
Throughout the present paper, we make this normalization. By the normalization and (\ref{Eigenvalue-1}), we see that $K'$ in (\ref{Eigenvalue-1}) must be zero. Hence, we obtain
\begin{equation}\label{Eigenvalue-2}
\Delta_{f} f + 2 \lambda f = 0.
\end{equation}

Next, by the first equation in (\ref{Ricci-Harmonic-Soliton}), we see that
\begin{equation}\label{Ineq-1}
\operatorname{Ric}_{g} + \operatorname{Hess} f = \lambda g + \alpha \nabla \phi \otimes \nabla \phi \geqslant \lambda g.
\end{equation}
We use the following theorem to obtain (\ref{Main-Theorem-Ineq}).

\begin{thm}[Futaki-Li-Li \cite{Futaki-Li-Li}]\label{Futaki-Li-Li-Theorem}
Let $(M, g)$ be a compact Riemannian manifold and $f : M \rightarrow \mathbb{R}$ a $\mathcal{C}^{2}$-class function on $M$. Suppose that
\[
\operatorname{Ric}_{g} + \operatorname{Hess} f \geqslant \lambda g
\]
for some real constant $\lambda \in \mathbb{R}$. Then the first non-zero eigenvalue $\lambda_{1}$ of the Witten-Laplacian $\Delta_{f}$ has the lower bound
\begin{equation}\label{Ineq-2}
\lambda_{1} \geqslant \sup_{s \in (0, 1)} \left \{ 4s(1 - s) \frac{\pi^{2}}{d^{2}} + s \lambda
\right \}, 
\end{equation}
where $d := \mathrm{diam}(M, g)$ denotes the diameter of $(M, g)$.
\end{thm}

The inequality (\ref{Ineq-1}) shows that Theorem \ref{Futaki-Li-Li-Theorem} above works for the compact domain manifold $M$ of the Ricci-harmonic soliton $((M, g), (N, h), \phi, f, \lambda)$. Hence, the first non-zero eigenvalue $\lambda_{1}$ of the Witten-Laplacian $\Delta_{f}$ has the lower bound (\ref{Ineq-2}). Recall from (\ref{Eigenvalue-2}) that $2 \lambda$ is an eigenvalue of the Witten-Laplacian $\Delta_{f}$. Hence, by (\ref{Eigenvalue-2}) and (\ref{Ineq-2}), for any $0 < s < 1$, we have
\[
2 \lambda \geqslant 4s(1 - s) \frac{\pi^{2}}{d^{2}} + s \lambda, 
\]
which yields
\[
\lambda \geqslant \frac{4s(1 - s)}{2 - s} \cdot \frac{\pi^{2}}{d^{2}}.
\]
An elementary calculation shows that
\[
\frac{4s(1 - s)}{2 - s} \leqslant 4(\sqrt{2} - 1)^{2}, 
\]
where the equality is attained for $s = 2 - \sqrt{2} \in (0, 1)$. Hence, we have
\[
\lambda \geqslant 4(\sqrt{2} - 1)^{2} \frac{\pi^{2}}{d^{2}}, 
\]
equivalently (\ref{Main-Theorem-Ineq}). The proof of Theorem \ref{Main-Theorem} is completed.
\end{proof}

\section{Concluding Remarks}

For any Ricci-harmonic soliton $((M, g), (N, h), \phi, f, \lambda)$ satisfying (\ref{Ricci-Harmonic-Soliton}), we denote by $S := R - \alpha | \nabla \phi |^{2}$ the trace of the symmetric two-tensor $\operatorname{Ric}_{g} - \alpha \nabla \phi \otimes \nabla \phi$. We put
\[
\begin{aligned}
c & := \inf_{x \in M} \{ \operatorname{Ric}_{g}(v, v) - \alpha \nabla \phi \otimes \nabla \phi(v, v) : v \in T_{x} M, | v | = 1 \}, \\
C & := \sup_{x \in M} \{ \operatorname{Ric}_{g}(v, v) - \alpha \nabla \phi \otimes \nabla \phi(v, v) : v \in T_{x} M, | v | = 1 \}.
\end{aligned}
\]

Fern\'{a}ndez-L\'{o}pez and Garc\'{i}a-R\'{i}o \cite{Fernandez-Lopez-Garcia-Rio} gave some lower diameter bounds for compact shrinking gradient Ricci solitons depending on the scalar and Ricci curvatures as well as on the range of the potential function. In Section \ref{Section-2}, by replacing the Ricci tensor $\operatorname{Ric}_{g}$ with the symmetric two-tensor $\operatorname{Ric}_{g} - \alpha \nabla \phi \otimes \nabla \phi$ in the argument of \cite{Futaki-Sano, Futaki-Li-Li}, we gave a lower diameter bound for compact domain manifolds of shrinking Ricci-harmonic solitons. By using this way, we easily see that the same lower diameter bounds as in \cite{Fernandez-Lopez-Garcia-Rio} also hold for compact domain manifolds of shrinking Ricci-harmonic solitons:

\begin{thm}
Let $((M, g), (N, h), \phi, f, \lambda)$ be a shrinking Ricci-harmonic soliton satisfying {\rm (\ref{Ricci-Harmonic-Soliton})}. Suppose that the domain manifold $M$ is compact. Then the diameter of the domain manifold $(M, g)$ has the lower bound
\[
\operatorname{diam}(M, g) \geqslant \max \left \{ \sqrt{\frac{2(f_{\mathrm{max}} - f_{\mathrm{min}})}{\lambda - c}}, \sqrt{\frac{2(f_{\mathrm{max}} - f_{\mathrm{min}})}{C - \lambda}}, 2 \sqrt{\frac{2(f_{\mathrm{max}} - f_{\mathrm{min}})}{C - c}} \right \}.
\]
\end{thm}

\begin{cor}
Let $((M, g), (N, h), \phi, f, \lambda)$ be a shrinking Ricci-harmonic soliton satisfying {\rm (\ref{Ricci-Harmonic-Soliton})}. Suppose that the domain manifold $M$ is compact. Then the diameter of the domain manifold $(M, g)$ has the lower bound
\[
\operatorname{diam}(M, g) \geqslant \max \left \{ \sqrt{\frac{S_{\mathrm{max}} - m \lambda}{\lambda(\lambda - c)}}, \sqrt{\frac{S_{\mathrm{max}} - m \lambda}{\lambda(C - \lambda)}}, 2 \sqrt{\frac{S_{\mathrm{max}} - m \lambda}{\lambda(C - c)}} \right \}.
\]
\end{cor}

\begin{cor}
Let $((M, g), (N, h), \phi, f, \lambda)$ be a shrinking Ricci-harmonic soliton satisfying {\rm (\ref{Ricci-Harmonic-Soliton})}. Suppose that the domain manifold $M$ is compact and $\operatorname{Ric}_{g} - \alpha \nabla \phi \otimes \nabla \phi > 0$. Then the diameter of the domain manifold $(M, g)$ has the lower bound
\[
\operatorname{diam}(M, g) \geqslant \max \left \{ \sqrt{\frac{S_{\mathrm{max}} - S_{\mathrm{min}}}{\lambda(\lambda - c)}}, \sqrt{\frac{S_{\mathrm{max}} - S_{\mathrm{min}}}{\lambda(C - \lambda)}}, 2 \sqrt{\frac{S_{\mathrm{max}} - S_{\mathrm{min}}}{\lambda(C - c)}} \right \}.
\]
\end{cor}

\end{document}